\documentclass[reqno]{amsart}
\usepackage{amssymb,hyperref}

\newtheorem{thm}{Theorem}
\newtheorem{lem}{Lemma}
\newtheorem{cor}{Corollary}

\numberwithin{equation}{section}

\renewcommand{\Re}{\mathrm{Re}}
\renewcommand{\Im}{\mathrm{Im}}

\title
[Starlikeness problems for certain analytic functions ...]
{Starlikeness problems \\
for certain analytic functions \\
concerned with subordinations}

\author
[H. Shiraishi, S. Owa, T. Hayami, K. Kuroki and H. M. Srivastava]
{Hitoshi Shiraishi, Shigeyoshi Owa, Toshio Hayami, \\
Kazuo Kuroki and H. M. Srivastava}

\address{Hitoshi Shiraishi \newline
Department of Mathematics \newline
Kinki University \newline
Higashi-Osaka, Osaka 577-8502, Japan}
\email{shiraishi@math.kindai.ac.jp}

\address{Shigeyoshi Owa \newline
Department of Mathematics \newline
Kinki University \newline
Higashi-Osaka, Osaka 577-8502, Japan}
\email{owa@math.kindai.ac.jp}

\address{Toshio Hayami \newline
School of Science and Technology \newline
Kinki University \newline
Sanda, Hyogo 669-1337, Japan}
\email{ha\_ya\_to112@hotmail.com}

\address{Kazuo Kuroki \newline
Department of Mathematics \newline
Kinki University \newline
Higashi-Osaka, Osaka 577-8502, Japan}
\email{freedom@sakai.zaq.ne.jp}

\address{H. M. Srivastava \newline
Departmemt of Mathematics and Statistics \newline
University of Victoria \newline
Victoria, British Columbia V8W 3P4, Canada}
\email{harimsri@math.uvic.ca}


\date{}

\begin{document}

\begin{abstract}
Let $\mathcal{A}_n$ be the class of functions $f(z)$ which are analytic in the open unit disk $\mathbb{U}$ with $f(0)=0$, $f'(0)=1$, $f''(0)=f'''(0)=\ldots=f^{(n)}=0$ and $f^{(n+1)}\neq0$.
Applying the results due to S. S. Miller (J. Math. Anal. Appl.{\bf 65}(1978), 289-305), some interesting starlikeness problems concerned with subordinations are discussed.
The results in the paper are extensions of results by M. Obradovi\'{c} (Hokkaido Math. J. {\bf 27}(1998), 329-335).
\end{abstract}

\begin{flushleft}
This paper was published in the journal: \\
Comput. Math. Appl. {\bf 62} (2011), No. 8, 2978--2987. \\
\url{http://www.sciencedirect.com/science/article/pii/S0898122111006675}
\end{flushleft}
\hrule

\

\

\maketitle

\section{Introduction}

\

Let $\mathcal{H}[a_0,n]$ denote the class of functions $p(z)$ of the form
$$
p(z)
= a_0 + \sum_{k=n}^{\infty} a_k z^k
\qquad (n=1,2,3,\ldots)
$$
which are analytic in the open unit disk $\mathbb{U}=\{z \in \mathbb{C}:|z|<1\}$,
where $a_0\in\mathbb{C}$.

Also,
let $\mathcal{A}_n$ denote the class of functions
$$
f(z)=z+a_{n+1}z^{n+1}+a_{n+2}z^{n+2}+ \ldots
\qquad(n=1,2,3,\ldots)
$$
that are analytic in $\mathbb{U}$ with $a_{n+1}\neq0$
and $\mathcal{A}\equiv\mathcal{A}_1$.

If $f(z)\in\mathcal{A}_n$ satisfies
$$
\Re \left( \frac{zf'(z)}{f(z)} \right)
> \alpha
\qquad (z\in\mathbb{U})
$$
for some real $\alpha$ $(0\leqq\alpha<1)$,
then we say that $f(z)$ is starlike of order $\alpha$ and written by $f(z)\in\mathcal{S}^*(\alpha)$
and $\mathcal{S}^*\equiv\mathcal{S}^*(0)$. 

Let $f(z)$ and $g(z)$ be analytic in $\mathbb{U}$.
Then $f(z)$ is said to be subordinate to $g(z)$ if there exists an analytic function $w(z)$ in $\mathbb{U}$ satisfying $w(0)=0$, $|w(z)| < 1$ $(z\in\mathbb{U})$ and such that $f(z)=g(w(z))$.
We denote this subordination by
$$
f(z)
\prec g(z)
\qquad (z\in\mathbb{U}).
$$

In particular, if $g(z)$ is univalent in $\mathbb{U}$, then the subordination
$$
f(z)
\prec g(z)
\qquad (z\in\mathbb{U})
$$
is equivalent to $f(0)=g(0)$ and $f(\mathbb{U}) \subset g(\mathbb{U})$ (cf. \cite{d7ref2}).

\

The basic tool in proving our results is the following lemma due to Miller and Mocanu \cite{m1ref2} (also \cite{d7ref2}).

\

\begin{lem} \label{jack} \quad
Let the function $w(z)$ defined by
$$
w(z)=a_nz^n+a_{n+1}z^{n+1}+a_{n+2}z^{n+2}+ \ldots
\qquad(n=1,2,3,\ldots)
$$
be analytic in $\mathbb{U}$ with $w(0)=0$.
If $\left|w(z)\right|$ attains its maximum value on the circle $\left|z\right|=r$ at a point $z_{0}\in\mathbb{U}$,
then there exists a real number $k \geqq n$ such that
$$
\frac{z_{0}w'(z_{0})}{w(z_{0})}=k.
$$
\end{lem}

\

\section{Main result}

\

Applying Lemma \ref{jack},
we have the following lemma.

\

\begin{lem} \label{d8lem1} \quad
Let $p(z)\in\mathcal{H}[1,n]$ satisfy the condition
$$
p(z)-\frac{1}{\mu}zp'(z)
\prec 1+\lambda z
\qquad (z\in\mathbb{U})
$$
for some complex number $\mu$ $(\mathrm{Re}(\mu)<n,\mu\neq0)$ and some complex number $\lambda$ $(0<|\lambda|\leqq1)$,
then
$$
p(z)
\prec 1+\lambda_1 z
\qquad (z\in\mathbb{U}),
$$
where $\lambda_1$ is a complex number such that
\begin{equation} \label{d8lem1eq1}
|\lambda_1|
= |\lambda| \frac{|\mu|}{|n-\mu|}.
\end{equation}
\end{lem}

\

\begin{proof} \quad
We consider the function $p(z)$ defined by
$$
p(z)
= 1+\lambda_1w(z)
$$
with $\lambda_1$ is given by (\ref{d8lem1eq1}).
Then, $w(z)$ is analytic in $\mathbb{U}$ and $w(0)=0$.
We want to show that $|w(z)|<1$ $(z\in\mathbb{U})$.
If there exists a point $z_0\in\mathbb{U}$ such that $|w(z_0)|=1$,
then we can write $z_0w'(z_0)=kw(z_0)$ $(k\geqq n)$ by Lemma \ref{jack}.
Putting $w(z_0)=e^{i\theta}$,
we get
\begin{align*}
\left| p(z_0)-\frac{1}{\mu}z_0p'(z_0)-1 \right|
&= \left| \lambda_1w(z_0)-\frac{1}{\mu}\lambda_1z_0w'(z_0) \right| \\
&= \left| \lambda_1e^{i\theta}-\frac{1}{\mu}\lambda_1ke^{i\theta} \right| \\
&= |\lambda_1|\frac{|\mu-k|}{|\mu|} \\
&= |\lambda_1|\frac{\sqrt[]{(k-\Re(\mu))^2+(\Im(\mu))^2}}{|\mu|} \\
&\geqq |\lambda_1|\frac{\sqrt[]{(n-\Re(\mu))^2+(\Im(\mu))^2}}{|\mu|} \\
&= |\lambda_1|\frac{|n-\mu|}{|\mu|}
= |\lambda|.
\end{align*}

This contradicts the assumption of Lemma \ref{d8lem1}.
Therefore, there is no $z_0\in\mathbb{U}$ such that $|w(z_0)|=1$.
This implies that $|w(z)|<1$ $(z\in\mathbb{U})$.
This completes the proof of the lemma.
\end{proof}

\

Next, we show

\

\begin{lem} \label{d8lem2} \quad
Let $\lambda$ and $\lambda_1$ be complex numbers such that $0<|\lambda_1|<|\lambda|<1$ and let $Q(z)\in\mathcal{H}[1,n]$ such that
\begin{equation} \label{d8lem2eq1}
Q(z)
\prec 1+\lambda_1z
\qquad (z\in\mathbb{U}).
\end{equation}
\begin{enumerate}
\item \label{d8lem2emu1} \qquad
If $p(z)\in\mathcal{H}[1,n]$ and
\begin{equation} \label{d8lem2eq2}
Q(z)[\alpha+(1-\alpha)p(z)]
\prec 1+\lambda z
\qquad (z\in\mathbb{U})
\end{equation}
for some real $\alpha$ such that
\begin{equation} \label{d8lem2eq3}
\alpha \leqq \left\{
\begin{array}{cl}
\dfrac{1-|\lambda|}{1+|\lambda_1|} & (0<|\lambda|+|\lambda_1|\leqq1) \\
\dfrac{1-(|\lambda|^2+|\lambda_1|^2)}{2(1-|\lambda_1|^2)} & (|\lambda|^2+|\lambda_1|^2\leqq1\leqq|\lambda|+|\lambda_1|),
\end{array}
\right.
\end{equation}
then $\Re(p(z))>0$.
\item \label{d8lem2emu2} \qquad
If $w(z)\in\mathcal{H}[0,n]$ and
\begin{equation} \label{d8lem2eq4}
Q(z)[1+w(z)]
\prec 1+\lambda z
\qquad (z\in\mathbb{U}),
\end{equation}
then
\begin{equation} \label{d8lem2eq5}
|w(z)|
< \frac{|\lambda|+|\lambda_1|}{1-|\lambda_1|}
\leqq 1,
\end{equation}
where $|\lambda|+2|\lambda_1|\leqq1$.
\end{enumerate}

The bound $(\ref{d8lem2eq5})$ and the value of $\alpha$ given by $(\ref{d8lem2eq3})$ are best possible.
\end{lem}

\

\begin{proof} \quad
To prove (\ref{d8lem2emu1}),
we think about
$$
g(z)
= \frac{1+\lambda z}{1+\lambda_1 e^{i\theta} z}
\qquad (z\in\mathbb{U})
$$
for some real $\theta$.

It follows that
$$
\left| g(z) - \frac{1-\lambda\overline{\lambda_1} r^2 e^{-i\theta}}{1-|\lambda_1|^2 r^2} \right|
\leqq \frac{|\lambda-\lambda_1 e^{i\theta}|r}{1-|\lambda_1|^2 r^2}
$$
for $|z|\leqq r<1$.

If we put $\lambda=|\lambda|e^{i\theta_0}$, $\lambda_1=|\lambda_1|e^{i\theta_1}$ and $\phi=\theta_0-\theta_1-\theta$,
we can rewrite
$$
\left| g(z) - \frac{1-|\lambda||\lambda_1|r^2 e^{i\phi}}{1-|\lambda_1|^2 r^2} \right|
\leqq \frac{||\lambda_1|-|\lambda| e^{i\phi}|r}{1-|\lambda_1|^2 r^2}
\qquad(|z| \leqq r < 1).
$$

Hence for all $r$ $(0<r<1)$,
we obtain
$$
\mathrm{Re}(g(z))
> \frac{1-|\lambda||\lambda_1|\cos\phi- \sqrt[]{|\lambda|^2-2|\lambda||\lambda_1|\cos\phi+|\lambda_1|^2}}{1-|\lambda_1|^2}
\equiv h(\phi).
$$

Note that
$$
h'(\phi) = \frac{|\lambda||\lambda_1|\sin\phi}{1-|\lambda_1|^2}\left(1 - \frac{1}{\sqrt{|\lambda|^2-2|\lambda||\lambda_1|\cos\phi +|\lambda_1|^2}}\right).
$$

Therefore, $h'(\phi) = 0$ if $\sin(\phi)=0$ or $\sqrt{|\lambda|^2-2|\lambda||\lambda_1|\cos\phi +|\lambda_1|^2}=1$. Since
$$
|\lambda|^2-2|\lambda||\lambda_1|\cos\phi +|\lambda_1|^2 \leqq \left(|\lambda| + |\lambda_1|\right)^2,
$$
if $0 < |\lambda|+|\lambda_1|\leqq 1$, then $h'(\phi) = 0$ for $\phi=2k\pi$ $(k=0, \pm 1, \pm 2, \cdots)$ and $\phi = (2k+1)\pi$ $(k=0, \pm 1, \pm 2, \cdots)$.
With this conditions,
we see that
$$
h(\phi) \geqq h(2k+1)\pi = \frac{1-|\lambda|}{1+|\lambda_1|} \geqq \alpha.
$$
If $\sqrt{|\lambda|^2-2|\lambda||\lambda_1|\cos\phi +|\lambda_1|^2}=1$,
then $h'(\phi)=0$ for $\phi = \phi_1$ such that
$$
|\lambda|^2 - |\lambda||\lambda_1|\cos\phi_1 +|\lambda_1|^2 = 1.
$$

Thus, we have that
$$
h(\phi) \geqq h(\phi_1) = \frac{1-(|\lambda|^2 + |\lambda_1|^2)}{2(1-|\lambda_1|^2)} \geqq \alpha
$$
with $|\lambda|^2 + |\lambda_1|^2 \leqq 1 \leqq |\lambda| + |\lambda_1|$.

On the other hand,
in view of (\ref{d8lem2eq1}) and (\ref{d8lem2eq2}),
we have
$$
\alpha+(1-\alpha)p(z)
\prec g(z)
\qquad (z\in\mathbb{U}).
$$

So,
we can lead $\mathrm{Re}(p(z))>0$.

To prove (\ref{d8lem2emu2}),
note that in view of (\ref{d8lem2eq1}) and (\ref{d8lem2eq4}) we have
\begin{align*}
|w(z)|
&= \left| \frac{(Q(z)(1+w(z))-1)-(Q(z)-1)}{Q(z)} \right| \\
&\leqq \frac{|Q(z)(1+w(z))-1|+|Q(z)-1|}{|Q(z)|} \\
&\leqq \frac{|\lambda|+|\lambda_1|}{|1+\lambda_1 z|} \\
&< \frac{|\lambda|+|\lambda_1|}{1-|\lambda_1|} \leqq 1
\qquad(z\in\mathbb{U}).
\end{align*}

Let us show the sharpness.

We take
$$
Q(z)
= 1+\lambda_1 w_1(z)
\qquad (z\in\mathbb{U})
$$
where $w_1(z)$ is analytic in $\mathbb{U}$ such that $w_1(0)=0$ and $|w_1(z)|<1$ $(z\in\mathbb{U})$.

In the case of (\ref{d8lem2emu1}),
we put
$$
Q(z)[\alpha+(1-\alpha)p(z)]
= 1+\lambda w_0(z)
\qquad (z\in\mathbb{U})
$$
where $w_0(z)$ is analytic in $\mathbb{U}$ such that $w_0(0)=0$ and $|w_0(z)|<1$ $(z\in\mathbb{U})$.
Since
$$
p(z) = \frac{1}{1-\alpha}\left(\frac{1+\lambda w_0(z)}{1+\lambda_1 w_1(z)} - \alpha\right),
$$
if we consider $w_0(z)$ such that $\lambda w_0(z)=-|\lambda|$ and $\lambda_1 w_1(z) = |\lambda_1|$, then we have
$$
\alpha = \frac{1-|\lambda|}{1+|\lambda_1|} \qquad(0 < |\lambda|+|\lambda_1|\leqq 1).
$$

Furthermore, if we take $w_0(z)$ such that $\lambda w_0(z) = |\lambda|e^{i\theta_0}$ and $w_1(z)$ such that $\lambda_1 w_(z)=|\lambda_1|e^{i\theta_1}$, then we obtain that
$$
\alpha = \frac{1-(|\lambda|^2 + |\lambda_1|^2)}{2(1-|\lambda_1|^2)}
$$
with
$$
|\lambda|\sin\theta_0 - |\lambda_1|\sin\theta_1 + |\lambda||\lambda_1|\sin(\theta_0-\theta_1)
=0
$$
and
\begin{align*}
& 2(1-|\lambda_1|^2)|\lambda|\cos\theta_0 + 2|\lambda|^2|\lambda_1|\cos\theta_1 + 2(1-|\lambda_1|^2)|\lambda||\lambda_1|\cos(\theta_0-\theta_1) \\
&= (|\lambda_1|^2-1)(1+|\lambda|^2-|\lambda_1|^2).
\end{align*}

In the case of (\ref{d8lem2emu2}),
we put
$$
Q(z)[1+w(z)]
= 1+\lambda w_0(z)
\qquad (z\in\mathbb{U})
$$
where $w_0(z)$ that $w_0(0)=0$ and $|w_0(z)|<1$ $(z\in\mathbb{U})$.

Then,
we can obtain
$$
w(z)
= \frac{\lambda w_0(z) - \lambda_1 w_1(z)}{1+\lambda w_0(z)}
\qquad (z\in\mathbb{U}).
$$

If we take
$$
w(z)
= \frac{|\lambda|+|\lambda_1|}{1-|\lambda_1|}
\qquad (z\in\mathbb{U}),
$$
then we have $\lambda w_0(z) = |\lambda|$ and $\lambda_1 w_1(z) = -|\lambda_1|$.
\end{proof}

\

If we consider some real $\lambda$, $\lambda_1$ and fixed $\alpha$ in Lemma \ref{d8lem2},
we obtain Corollary \ref{d8cor4} due to S. Ponnusamy and V. Singh \cite{d8ref5}.

\

\begin{cor} \label{d8cor4} \quad
Let $\lambda$ and $\lambda_1$ be real with $0<\lambda_1<\lambda<1$ and let $Q(z)\in\mathcal{H}[1,n]$ satisfy
$$
Q(z)
\prec 1+\lambda_1z
\qquad (z\in\mathbb{U}).
$$
\begin{enumerate}
\item \qquad
If $p(z)\in\mathcal{H}[1,n]$ and
$$
Q(z)[\alpha+(1-\alpha)p(z)]
\prec 1+\lambda z
\qquad (z\in\mathbb{U})
$$
for some real $\alpha$ such that
\begin{equation} \label{d8cor4eq1}
\alpha = \left\{
\begin{array}{cl}
\dfrac{1-\lambda}{1+\lambda_1} & (0<\lambda+\lambda_1\leqq1) \\
\dfrac{1-(\lambda^2+\lambda_1^2)}{2(1-\lambda_1^2)} & (\lambda^2+\lambda_1^2\leqq1\leqq\lambda+\lambda_1),
\end{array}
\right.
\end{equation}
then $\Re(p(z))>0$.
\item \qquad
If $w(z)\in\mathcal{H}[0,n]$ and
$$
Q(z)[1+w(z)]
\prec 1+\lambda z
\qquad (z\in\mathbb{U}),
$$
then
\begin{equation} \label{d8cor4eq2}
|w(z)|
< \frac{\lambda+\lambda_1}{1-\lambda_1}
\leqq 1,
\end{equation}
where $\lambda+2\lambda_1\leqq1$.
\end{enumerate}

The bound $(\ref{d8cor4eq2})$ and the value of $\alpha$ given by $(\ref{d8cor4eq1})$ are best possible.
\end{cor}

\

By virtue of Lemma \ref{d8lem1},
we deduce the sufficient condition for the class $\mathcal{S}^*$.

\

\begin{thm} \label{d8thm1} \quad
If $f(z)\in\mathcal{A}_n$ satisfies the condition
$$
f'(z) \left( \frac{z}{f(z)} \right)^{1+\mu}
\prec 1+\lambda z
\qquad (z\in\mathbb{U})
$$
for some complex numbers $\mu$ $(\Re(\mu)<n)$ and $\lambda$ such that \\
$0<|\lambda|\leqq\dfrac{|n-\mu|}{\sqrt[]{|n-\mu|^2+|\mu|^2}}$,
then $f(z)\in\mathcal{S}^*$.
\end{thm}

\

\begin{proof} \quad
If we put
$$
Q(z)
= \left( \frac{z}{f(z)} \right)^\mu,
$$
then we get
$$
Q(z)-\frac{1}{\mu}Q'(z)
= f'(z)\left( \frac{z}{f(z)} \right)^{1+\mu}
\prec 1+\lambda z
\qquad (z\in\mathbb{U})
$$
for $0<|\lambda|\leqq1$.

In view of Lemma \ref{d8lem1}, we obtain
\begin{equation} \label{d8thm1eq1}
Q(z)
\prec 1+\lambda_1 z
\qquad (z\in\mathbb{U})
\end{equation}
and
$$
|\lambda_1|
= |\lambda| \frac{|\mu|}{|n-\mu|}.
$$
Then, we see that
$$
\left| \arg \left( f'(z)\left( \frac{z}{f(z)} \right)^{1+\mu} \right) \right|
< \arctan \left( \frac{|\lambda|}{\sqrt[]{1-|\lambda|^2}} \right)
$$
and
$$
\left| \arg \left( \left( \frac{f(z)}{z} \right)^\mu \right) \right|
= \left| \arg \left( \left( \frac{z}{f(z)} \right)^\mu \right) \right|
< \arctan \left( \frac{|\lambda_1|}{\sqrt[]{1-|\lambda_1|^2}} \right),
$$
which give
\begin{align*}
\left| \arg \left( \frac{zf'(z)}{f(z)} \right) \right|
&\leqq \left| \arg \left( f'(z)\left( \frac{z}{f(z)} \right)^{1+\mu} \right) \right| + \left| \arg \left( \left( \frac{f(z)}{z} \right)^\mu \right) \right| \\
&< \arctan \left( \frac{|\lambda|}{\sqrt[]{1-|\lambda|^2}} \right) + \arctan \left( \frac{|\lambda_1|}{\sqrt[]{1-|\lambda_1|^2}} \right) \\
&= \arctan \left( \dfrac{\dfrac{|\lambda|}{\sqrt[]{1-|\lambda|^2}} + \dfrac{|\lambda_1|}{\sqrt[]{1-|\lambda_1|^2}}}{1 - \dfrac{|\lambda||\lambda_1|}{\sqrt[]{(1-|\lambda|^2)(1-|\lambda_1|^2)}}} \right)
\leqq \frac{\pi}{2}
\qquad (z\in\mathbb{U}).
\end{align*}

This implies that $f(z)\in\mathcal{S^*}$

On the other hand,
when we have
$$
1 - \dfrac{|\lambda||\lambda_1|}{\sqrt[]{(1-|\lambda|^2)(1-|\lambda_1|^2)}}
\geqq 0,
$$
$\lambda$ satisfies
$$
|\lambda|
\leqq \frac{|n-\mu|}{\sqrt[]{|n-\mu|^2+|\mu|^2}}
< 1.
$$
\end{proof}

\

Taking $n=1$, $0<\mu<1$ and $0<\lambda<1$,
we have the next corollary due to Obradovi\'{c} \cite{d8ref4}.

\

\begin{cor} \label{d8cor1} \quad
If $f(z)\in\mathcal{A}$ satisfies the condition
$$
\left| f'(z) \left( \frac{z}{f(z)} \right)^{1+\mu}-1 \right|
< \lambda
\qquad (z\in\mathbb{U})
$$
for some $0<\mu<1$ and $0<\lambda<1$,
then $f(z)\in\mathcal{S}^*$.
\end{cor}

\

Applying Lemma \ref{d8lem2},
we derive the following theorem.

\

\begin{thm} \label{d8thm2} \quad
Let $f(z)\in\mathcal{A}_n$ satisfy the condition
\begin{equation} \label{d8thm2eq1}
f'(z) \left( \frac{z}{f(z)} \right)^{1+\mu}
\prec 1+\lambda z
\qquad (z\in\mathbb{U})
\end{equation}
for some complex number $\mu$ $\left( \Re(\mu)<\dfrac{n}{2} \right)$.
If the complex number $\lambda_1$ is given by
$$
|\lambda_1|
= |\lambda| \frac{|\mu|}{|n-\mu|},
$$
then
\begin{enumerate}
\item \qquad
$f(z)\in\mathcal{S}^*(\alpha)$ where
$$
\alpha \leqq \left\{
\begin{array}{ll}
\dfrac{1-|\lambda|}{1+|\lambda_1|}
& \left( 0<|\lambda|\leqq\dfrac{|n-\mu|}{|n-\mu|+|\mu|} \right) \\
\dfrac{1-(|\lambda|^2+|\lambda_1|^2)}{2(1-|\lambda_1|^2)}
& \left( \dfrac{|n-\mu|}{|n-\mu|+|\mu|}\leqq|\lambda|\leqq\dfrac{|n-\mu|}{\sqrt[]{|n-\mu|^2+|\mu|^2}} \right)
\end{array}
\right.
$$
and
\item \qquad
$\left| \dfrac{zf'(z)}{f(z)}-1 \right|<\dfrac{(|n-\mu|+|\mu|)|\lambda|}{|n-\mu|-|\mu||\lambda|}\leqq1$
\quad $(z\in\mathbb{U})$, \\
where $0<|\lambda|\leqq\dfrac{|n-\mu|}{|n-\mu|+2|\mu|}$.
\end{enumerate}
\end{thm}

\

\begin{proof} \quad
Let us define
$$
Q(z)
= \left( \frac{z}{f(z)} \right)^\mu,
$$
$$
p(z)
= \frac{zf'(z)}{f(z)},
$$
and
$$
w(z)
= \frac{zf'(z)}{f(z)}-1.
$$

Then by the equation (\ref{d8thm1eq1}) of Theorem \ref{d8thm1},
we have
$$
Q(z)
\prec 1+\lambda_1z
\qquad (z\in\mathbb{U})
$$
where $0 < |\lambda_1| = |\lambda|\dfrac{|\mu|}{|n-\mu|} < |\lambda| < 1$
since $\mathrm{Re}(\mu)<\dfrac{n}{2}$.
Also, since the condition (\ref{d8thm2eq1}) is equivalent to
$$
Q(z) \left[ \alpha+(1-\alpha)\frac{p(z)-\alpha}{1-\alpha} \right]
\prec 1+\lambda z
\qquad (z\in\mathbb{U})
$$
with $\alpha$ is given by (\ref{d8lem2eq3}) and as
$$
Q(z)[1+w(z)]
\prec 1+\lambda z
\qquad (z\in\mathbb{U}),
$$
the statement of the theorem directly follows from Lemma \ref{d8lem2}.
\end{proof}

\

Putting $n=1$, $0<\mu<\dfrac{1}{2}$, $0<\lambda<1$ and fixed $\alpha$,
we get the following corollary due to Obradovi\'{c} \cite{d8ref4}.

\

\begin{cor} \label{d8cor2} \quad
Let $f(z)\in\mathcal{A}$ satisfy the condition
$$
\left| f'(z) \left( \frac{z}{f(z)} \right)^{1+\mu}-1 \right|
< \lambda
\qquad (z\in\mathbb{U})
$$
with some $0<\mu<\dfrac{1}{2}$.
If the real number $\lambda_1$ is given by
$$
\lambda_1
= \lambda \frac{\mu}{1-\mu},
$$
then
\begin{enumerate}
\item \qquad
$f(z)\in\mathcal{S}^*(\alpha)$, where
$$
\alpha = \left\{
\begin{array}{ll}
\dfrac{1-\lambda}{1+\lambda_1}
& (0<\lambda\leqq1-\mu) \\
\dfrac{1-(\lambda^2+\lambda_1^2)}{2(1-\lambda_1^2)}
& \left( 1-\mu\leqq\lambda\leqq\dfrac{1-\mu}{\sqrt[]{(1-\mu)^2+\mu^2}} \right).
\end{array}
\right.
$$
\item \qquad
$\left| \dfrac{zf'(z)}{f(z)}-1 \right|<\dfrac{\lambda}{1-\mu-\mu\lambda}\leqq1$
\quad $(z\in\mathbb{U})$ ,\\
where $0<\lambda\leqq\dfrac{1-\mu}{1+\mu}$.
\end{enumerate}
\end{cor}

\

Using Lemma \ref{d8lem1} and Lemma \ref{d8lem2},
we derive

\

\begin{thm} \label{d8thm3} \quad
If $f(z)\in\mathcal{A}_n$ satisfies
$$
f'(z) \left( \frac{z}{f(z)} \right)^{1+\mu}
\prec 1+ \lambda z
\qquad (z\in\mathbb{U})
$$
and
$$
F(z)
= z \left[ \frac{c-\mu}{z^{c-\mu}} \int_0^z \left( \frac{t}{f(t)} \right)^\mu t^{c-\mu-1} dt \right]^{-\frac{1}{\mu}}
\qquad (z\in\mathbb{U})
$$
for some complex numbers $\lambda, \mu$, and $c$ such that $\Re(c-\mu)<n$,
then
\begin{enumerate}
\item \qquad
$F(z)\in\mathcal{S}^*$ for $|c-\mu||\lambda|\leqq\dfrac{|n-\mu||n-(c-\mu)|}{\sqrt[]{|n-\mu|^2+|\mu|^2}}$.
\item \qquad
$F(z)\in\mathcal{S}^*(\alpha)$ where
$$
\alpha \leqq \left\{
\begin{array}{ll}
\dfrac{1-|\lambda_1|}{1+|\lambda_2|}
& \left( 0<|\lambda_1|\leqq\dfrac{|n-\mu|}{|n-\mu|+|\mu|} \right) \\
\dfrac{1-(|\lambda_1|^2+|\lambda_2|^2)}{2(1-|\lambda_2|^2)}
& \left( \dfrac{|n-\mu|}{|n-\mu|+|\mu|}\leqq|\lambda_1|\leqq\dfrac{|n-\mu|}{\sqrt[]{|n-\mu|^2+|\mu|^2}} \right),
\end{array}
\right.
$$
$|\lambda_1|=|\lambda|\dfrac{|c-\mu|}{|n-(c-\mu)|}$,
$|\lambda_2|=|\lambda_1|\dfrac{|\mu|}{|n-\mu|}$ and $\Re(\mu)<\dfrac{n}{2}$.
\item \qquad
$\left| \dfrac{zF'(z)}{F(z)}-1 \right|<\dfrac{(|n-\mu|+|\mu|)|c-\mu||\lambda|}{|n-(c-\mu)||n-\mu|-|c-\mu||\mu||\lambda|}\leqq1$
\quad $(z\in\mathbb{U})$
where $|c-\mu||\lambda|\leqq\dfrac{|n-(c-\mu)||n-\mu|}{|n-\mu|+2|\mu|}$
and $\Re(\mu)<\dfrac{n}{2}$.
\end{enumerate}
\end{thm}

\

\begin{proof} \quad
If we put
$$
Q(z)
= F'(z) \left( \frac{z}{F(z)} \right)^{1+\mu},
$$
after some transformations, we obtain
$$
Q(z)+\frac{1}{c-\mu}zQ'(z)
= f'(z)\left( \frac{z}{f(z)} \right)^{1+\mu}
\prec 1+\lambda z.
$$

Therefore, spending the same technique as in the proof of Lemma \ref{d8lem1},
we have that
$$
Q(z)
\prec 1+\lambda_2 z
\qquad (z\in\mathbb{U})
$$
and
$$
|\lambda_2|
= \frac{|\mu||c-\mu|}{|n-\mu||n-(c-\mu)|}.
$$

The statement of the theorem now easily follows from Theorem \ref{d8thm1} and Theorem \ref{d8thm2}.
\end{proof}

\

Taking $n=1$, $0<\mu<1$, $0<\lambda<1$ and fixed $\alpha$,
we obtain the next the corollary due to Obradovi\'{c} \cite{d8ref4}.

\

\begin{cor} \label{d8cor3} \quad
If $f(z)\in\mathcal{A}$ satisfies
$$
\left| f'(z) \left( \frac{z}{f(z)} \right)^{1+\mu}-1 \right|
< \lambda
\qquad (z\in\mathbb{U})
$$
and
$$
F(z)
= z \left[ \frac{c-\mu}{z^{c-\mu}} \int_0^z \left( \frac{t}{f(t)} \right)^\mu t^{c-\mu-1} dt \right]^{-\frac{1}{\mu}}
\qquad (z\in\mathbb{U})
$$
for $c-\mu>0$,
then
\begin{enumerate}
\item \qquad
$F(z)\in\mathcal{S}^*$ for $(c-\mu)\lambda\leqq\dfrac{(1-\mu)(1-(c-\mu))}{\sqrt[]{(1-\mu)^2+\mu^2}}$.
\item \qquad
$F(z)\in\mathcal{S}^*(\alpha)$ where
$$
\alpha = \left\{
\begin{array}{ll}
\dfrac{1-\lambda_1}{1+\lambda_2}
& (0<\lambda_1\leqq1-\mu) \\
\dfrac{1-(\lambda_1^2+\lambda_2^2)}{2(1-\lambda_2^2)}
& \left( 1-\mu\leqq\lambda_1\leqq\dfrac{1-\mu}{\sqrt[]{(1-\mu)^2+\mu^2}} \right),
\end{array}
\right.
$$
$\lambda_1=\lambda\dfrac{c-\mu}{1-(c-\mu)}$,
$\lambda_2=\lambda_1\dfrac{\mu}{1-\mu}$ and $0<\mu<\dfrac{1}{2}$.
\item \qquad
$\left| \dfrac{zF'(z)}{F(z)}-1 \right|<\dfrac{(c-\mu)\lambda}{(1-(c-\mu))(1-\mu)-(c-\mu)\mu\lambda}\leqq1$
\quad $(z\in\mathbb{U})$
where $0<\lambda\leqq\dfrac{(1-(c-\mu))(1-\mu)}{(c-\mu)(1+\mu)}$
and $0<\mu<\dfrac{1}{2}$.
\end{enumerate}
\end{cor}

\

\end{document}